\documentclass[12pt]{amsart}
\usepackage{amsmath, amssymb}

\newcommand{\PP}{\ensuremath{\mathbb{P}}}
\newcommand{\ZZ}{\ensuremath{\mathbb{Z}}}

\newcommand{\CC}{\ensuremath{{\mathbb{C}}}}
\newcommand{\QQ}{\ensuremath{{\mathbb{Q}}}}
\newcommand{\CS}{\ensuremath{\CC^{{*}}}}
\newcommand{\RR}{\ensuremath{\mathbb{R}}}

\DeclareMathOperator{\Sec}{Sec}

\DeclareMathOperator{\conv}{conv}

\newtheorem{thm}{Theorem}[section]
\newtheorem{cor}[thm]{Corollary}
\newtheorem{lem}[thm]{Lemma}
\newtheorem{prop}[thm]{Proposition}
\theoremstyle{definition}

\newtheorem{defin}[thm]{Definition}
\newtheorem{ex}[thm]{Example}
\newtheorem{rmk}[thm]{Remark}

\begin{document}

\title{Secant varieties of toric varieties}

\author{David Cox}
\address{Department of Mathematics and Computer Science, Amherst
  College, Amherst, MA 01002}
\email{dac@cs.amherst.edu}

\author{Jessica Sidman}
\address{Department of Mathematics and Statistics, Mount Holyoke
  College, South Hadley, MA 01075}
\email{jsidman@mtholyoke.edu}

\begin{abstract}
Let $X_P$ be a smooth projective toric variety of dimension $n$
embedded in $\PP^r$ using all of the lattice points of the polytope
$P$.  We compute the dimension and degree of the secant variety $\Sec
X_P$.  We also give explicit formulas in dimensions $2$ and $3$ and
obtain partial results for the projective varieties $X_A$ embedded
using a set of lattice points $A \subset P\cap\ZZ^n$ containing the
vertices of $P$ and their nearest neighbors.
\end{abstract}

\maketitle

\section{Introduction}
Let $X \subseteq \PP^r$ be a reduced and irreducible complex variety of 
dimension $n.$  Its
$k$th secant variety, $\Sec_k X,$ is the closure of the union of all
$(k-1)$-planes in $\PP^r$ meeting $X$ in at least $k$ points.
In this paper, we discuss the dimension and degree of the secant
variety $\Sec X = \Sec_2 X$ when $X$ is a smooth toric variety equivariantly
embedded in projective space.

Secant varieties arise naturally in classical examples (see 
Examples~\ref{ex:veronese}, \ref{ex:segre} and~\ref{ex:scroll} below) and 
also in the newly developing fields of
algebraic statististics and phylogenetic combinatorics (see, for
example, \cite{erss}).   The basic theory
of secant varieties is explained in the books \cite{fov,zak}.
 
There is a rich
collection of ideas relating secant varieties, tangent varieties, dual
varieties, and the Gauss map, as discussed in \cite{gh} and
\cite[Section 4.4]{fov}.  Many authors have studied the problem of classifying 
the varieties $X$ such that
$\Sec_k X$ has the expected dimension, $\min\{r, k(n+1) - 1\},$ and what
degeneracies may occur \cite{cj, cj2, cgg2, cgg, ohno}.  For the duals of toric varieties,
this was done in \cite{DiRocco}.  Our paper can be regarded as the
beginnings of a similar study for the secant varieties of toric
varieties.

Currently, there is much activity focused on refining our understanding 
of secant varieties.  A lower bound for the degree of $\Sec_k X$ is 
given in \cite{cr}, and the degree of the secant variety of a monomial curve
is worked out in \cite{r}.  Questions on the defining equations of secant 
varieties are discussed in \cite{lm,vermeire}.  The recent paper \cite{ss}
uses combinatorial Gr\"obner methods to study the ideal of $\Sec_k X$ and its
degree and dimension.

Many classical varieties whose secant varieties have been studied in
the literature are toric---see for example \cite{ballico, cj,
cj2,cgg,depoi, geramita,landsberg}.  The secant varieties below will
play an important role in our paper.

\begin{ex}
\label{ex:veronese}
If $X$ is the image of the Veronese map $\nu_2:\PP^n \to \PP^r$ by
$\mathcal{O}_{\PP^n}(2)$, then its secant variety has the expected
dimension only when $n = 1$; otherwise the dimension is $2n < \min\{r,
2n+1\}$.

The ideal of $X$ is generated by the $2 \times 2$
minors of the $(n+1) \times (n+1)$ generic symmetric matrix and
the ideal of $\Sec X$ is generated by the $3 \times 3$ minors of the
same matrix.  (See Example 1.3.6 in \cite{fov}.)
\end{ex}

\begin{ex}
\label{ex:segre}
If $X = \PP^\ell \times \PP^{n-\ell}$, $1 \le \ell \le n-1$, is
embedded in $\PP^r$ via the Segre map, then $\Sec X$ has the expected
dimension only when $\ell = 1,n-1$; otherwise $\dim \Sec X = 2n-1 <
\min\{r,2n+1\}$.

Here, the ideal of $X$ is also determinantal, generated by the $2
\times 2$ minors of the generic $(\ell+1) \times (n-\ell+1)$ matrix,
and the ideal of $\Sec X$ is defined by the $3 \times 3$ minors of the
same matrix when $\ell+1,n-\ell+1 \ge 3$.  (See Example 4.5.21 of
\cite{fov}.)
\end{ex}

\begin{ex}
\label{ex:scroll}
Given positive integers $d_1,\dots,d_n$, the rational normal scroll
$S_{d_1,\dots,d_n}$ is the image of
\[
\PP(\mathcal{E}),\quad \mathcal{E} = \mathcal{O}_{\PP^1}(d_1)\oplus
\dots \oplus \mathcal{O}_{\PP^1}(d_n)
\]
under the projective embedding given the ample line bundle
$\mathcal{O}_{\PP(\mathcal{E})}(1)$.  The variety $\Sec X$ always has the
expected dimension, which is $2n+1$ except when $\sum_{i=1}^n d_i \le
n+1$.

Furthermore, the ideal of $X$ (resp.\ $\Sec X$) is generated by the $2
\times 2$ (resp.\ $3 \times 3$) minors of a matrix built out of Hankel
matrices.  (See \cite{cj} and Section 4.7 of \cite{fov}.)
\end{ex}

A key observation is that Examples~\ref{ex:veronese}, \ref{ex:segre}
and~\ref{ex:scroll} each involve a smooth toric variety coming from a
particularly simple polytope.  It is natural ask what happens for more
general toric varieties.  To state our main theorem, we introduce
some definitions.

If $P \subset \RR^n$ is an $m$-dimensional lattice polytope, we say
that $P$ is \emph{smooth} if it is simplicial and the first lattice
vectors along the edges incident at any vertex form part of a
$\ZZ$-basis of $\ZZ^n$.  (Note that such polytopes are also called
\emph{Delzant} in the literature.)  The lattice points $\{u_0,\ldots,
u_r \} = P\cap\ZZ^n$ give characters $\chi^{u_i}$ of the torus
$(\CS)^n$ that give an embedding $X \hookrightarrow \PP^r$ defined by
$x \mapsto [\chi^{u_0}(x): \cdots : \chi^{u_r}(x)]$, where $X$ is the
abstract toric variety associated to the inner normal fan of $P$.  We
denote the image of $X$ under this embedding by $X_P$.  Note that $X$
is smooth if and only if $P$ is smooth.

We also observe that $X_P$ is projectively equivalent to $X_Q$ if
$P$ is obtained from $Q$ via an element of the group $AGL_n(\ZZ)$ of
affine linear isomorphisms of $\ZZ^n$.

The \emph{standard simplex} of dimension $n$ in $\RR^n$ is
\[
\Delta_n =  \{(p_1, \ldots, p_n) \in \RR^n_{\ge 0}
\mid p_1 + \cdots + p_n \le 1 \}
\]
and its multiple by $r \ge 0$ is
\[
r\Delta_n =  \{(p_1, \ldots, p_n) \in \RR^n_{\ge 0}
\mid p_1 + \cdots + p_n \le r \}.
\]

We let $(2\Delta_n)_k$ denote the convex hull of the lattice points in
$2\Delta_n$ minus a $k$-dimensional face, and let $Bl_k(\PP^n)$
denote the blowup of $\PP^n$ along a torus-invariant $k$-dimensional
linear subspace.

Here is the main result of the paper.

\begin{thm}\label{thm:dimension}
Let $P$ be a smooth polytope of dimension $n$.  The dimension and
degree of $\Sec X_P$ are
given in the following table:

\begin{table}[h]
\label{table:dimdeg}
{\small
\begin{equation*}
\begin{array}{|c|c|@{\ \ }l|}\hline
\hbox{\vrule height13pt depth3.5pt width0pt}P & X_P & \text{\rm
Dimension and degree of}\ \Sec X_P\\[2pt] 
\hline\hline
\hbox{\vrule height13pt depth3.5pt width0pt}\Delta_n & \PP^n & \dim
\Sec X_P = n\\[3pt] 
&& \deg \Sec X_P\, = 1\\[2pt]
\hline
\hbox{\vrule height13pt depth3.5pt width0pt} 2\Delta_n & \nu_2(\PP^n)
& \dim \Sec X_P = 2n\\[3pt] 
&&  \deg \Sec X_P\, = \binom{2n-1}{n-1} \\[4pt]
\hline
\hbox{\vrule height13pt depth3.5pt width0pt}(2\Delta_n)_k &
Bl_k(\PP^n) & \dim \Sec X_P = 2n \\[3.5pt] 
{\text{\footnotesize $(0 \le k \le n-2)$}} && \deg \Sec X_P\, = \\[1.5pt]
&& \quad\quad\quad{\displaystyle
\!\!\sum_{1\le i < j \le n-k}}\!\!\!\!
\bigl[\binom{n}{n-i}\binom{n-1}{n-j} -
\binom{n}{n-j}\binom{n-1}{n-i}\bigr]\\[14.5pt]  
\hline
\hbox{\vrule height13pt depth3.5pt width0pt} \Delta_{\ell} \times
\Delta_{n-\ell} & \PP^\ell \times \PP^{n-\ell} & \dim \Sec X_P = 2n-1\\[-.5pt]
\text{\footnotesize $(1 \le \ell \le n-1)$} & &  \deg
\Sec X_P\, =  
{\displaystyle
\prod_{0 \le i \le
n-\ell-2}}\frac{\binom{\ell+1+i}{2}}{\binom{2+i}{2}}\\[14pt]  
\hline
\hbox{\vrule height13pt depth3.5pt width0pt} \text{\rm not}\,
AGL_n(\ZZ)\text{-} & 
X_P & \dim \Sec X_P = 2n+1\\[3pt]
\text{\rm equivalent to} && \deg \Sec X_P\, = \\[-1pt]
\raise3.5pt\hbox{\text{\rm the above}} && {\text{\footnotesize
 $\,\frac12\big((\deg X_P)^2 - 
{\displaystyle\sum_{i = 0}^{n}}\binom{2n+1}{i}
{\displaystyle\int_{X_P}}\!\! c(T_{X_P})^{-1} \cap c_1(\mathcal{L})^i\big)
$}}\\[11pt]
\hline
\end{array}
\end{equation*}}
\end{table}

\noindent In the last row, $\mathcal{L} = \mathcal{O}_{X_P}(1)$ and
$c(T_{X_P})^{-1}$ is the inverse of the total Chern class of the
tangent bundle of $X_P$ in the Chow ring of $X_P$.  Furthermore:
\begin{enumerate}
\item For $P$ in the first four rows of the table, there are
infinitely many secant lines of $X_P$ through a general point of $\Sec
X_P$.
\item For $P$ in the last row of the table,
there is a unique secant line of $X_P$ through a general point of
$\Sec X_P$.
\end{enumerate}
\end{thm}

For the remainder of the paper, the table appearing in
Theorem~\ref{thm:dimension} will be referred to as
Table~\ref{table:dimdeg}.  In Section~\ref{sec:complete} we prove that
the first four rows of Table~\ref{table:dimdeg} give all smooth
subpolytopes of $2\Delta_n$ up to $AGL_n(\ZZ)$-equi\-valence.  Hence
$\dim \Sec X_P = 2n+1$ if and only if $P$ does not fit inside
$2\Delta_n$.  This also determines the number of secant lines through
a general point of $\Sec X_P$.

Theorem~\ref{thm:dimension} enables us to decide when the secant
variety has the expected dimension.

\begin{cor}\label{cor:dimension}
When $X_P \subset \PP^r$ comes from a smooth polytope $P$ of dimension
$n$, $\Sec X_P$ has the expected dimension $\min\{r,2n+1\}$ unless $P$
is $AGL_n(\ZZ)$-equivalent to one of
\[
2\Delta_n\ (n\ge2),\ (2\Delta_n)_k\ (0 \le k \le n-3),\
\Delta_\ell \times \Delta_{n-\ell}\ (2 \le \ell \le n-2).
\]
\end{cor}

We can compute the degree of $\Sec X_P$ explicitly in low dimensions
as follows.  Throughout, $n$-dimensional volume is normalized so that
the volume of $\Delta_n$ is $1$.

\begin{cor}\label{thm:sfdegree}
Let $P$ be a smooth polygon in $\RR^2$.  If $\dim \Sec X_P =
5$, then
\[ 
\deg \Sec X_P = \tfrac{1}{2}(d^2-10d+5B+2V-12), 
\]
where $d$ is the area of $P$, $B$ is the number of lattice points on
the boundary of $P$, and $V$ is the number of vertices of $P$.
\end{cor}

\begin{cor}\label{thm:3folddegree}
Let $P$ be a smooth 3-dimensional polytope in $\RR^3$.
If $\dim \Sec X_P = 7$, then
\[ 
\deg \Sec X_P = \tfrac{1}{2}(d^2-21d+c_1^3+8V+14E-84I-132), 
\]
where $d$ is the volume of $P,$ $E$ is the number of lattice points on
the edges of $P$, $V$ is the number of vertices of $P$, and $I$ is the
number of interior lattice points.  Also, $c_1 = c_1(T_{X_P})$.
\end{cor}

\begin{rmk}
Computing explicit degree formulas for $\Sec X_P$ in terms of the
lattice points of $P$ gets more difficult as the dimension increases.
However, we would like to stress that given any smooth polytope $P$,
the formula in the last row of {Table~\ref{table:dimdeg}} is
computable solely from the data of $P$.  Indeed, $c(T_{X_P})^{-1}$
depends only on the inner normal fan of $P$, and $\mathcal{L}$ is the
line bundle given by $P$.  The intersection products may be determined
from the combinatorial geometry of $P$ using standard results in
\cite{fulton-t} or any convenient presentation of the Chow ring of
$X_P$.  We give an example of this in Theorem~\ref{thm:segver} when we
compute the degree of the secant variety of a Segre-Veronese variety.
\end{rmk}

We prove Theorem~\ref{thm:dimension} in three steps as follows.

First, for the polytopes in the first four rows of
{Table~\ref{table:dimdeg}}, the degrees and dimensions of their secant
varieties are computed in Section~\ref{sec:subpolytopes} using
determinantal presentations for the ideals of $X_P$ and $\Sec X_P$,
together with known results about dimensions and degrees of
determinantal varieties.  We also relate these polytopes to
subpolytopes of $2\Delta_n$.

Second, the bottom row of {Table~\ref{table:dimdeg}} uses a well-known
formula for the degree of $\Sec X_P$ times the degree of the linear
projection from the abstract join of $X_P$ with itself to $\PP^r$.  We
discuss this formula in Section~\ref{sec:degformula} and apply it to
various examples.

Third, for the polytopes in the bottom row of
{Table~\ref{table:dimdeg}}, we prove in Section~\ref{sec:complete}
that a general point of $\Sec X_P$ lies on a unique secant line by
showing that $P$ contains configurations of lattice points that are
easier to study.  Known results about rational normal scrolls will be
used in the proof.

{}From these results, Theorem~\ref{thm:dimension} and
Corollaries~\ref{cor:dimension}, \ref{thm:sfdegree} and
\ref{thm:3folddegree} follow easily.  The paper concludes with
Section~\ref{sec:projections}, where we study the toric varieties
$X_A$ embedded using a subset of lattice points $A \subset P\cap\ZZ^n$
of a smooth polytope $P$ of dimension $n$.  When $A$ contains the
vertices of $P$ and their nearest neighbors along the edges, we
compute $\dim \Sec X_A$ and obtain partial results about $\deg \Sec
X_A$.

\subsection*{Acknowledgements} We are grateful to Bernd Sturmfels for
suggesting the problem and for helpful conversations.  We are also
grateful to Sturmfels and Seth Sullivant for pointing out the
importance of using vertices and their nearest neighbors along the
edges.  We thank Masahiro Ohno for a useful reference and anonymous
referees for several helpful suggestions and for prompting us to
strengthen Theorem~\ref{thm:dimension}.  We also thank Aldo Conca
for communications on determinantal ideals.  We made numerous experiments
using \emph{Macaulay 2} \cite{M2}, {\tt polymake} \cite{polymake},
Maple, and Mathematica, along with {\scshape Singular} \cite{GPS01}
code provided by Jason Morton.  The second author was supported by an
NSF postdoctoral fellowship, Grant No.\ 0201607, and also thanks the
University of Massachusetts, Amherst for their hospitality in
2004-05 and the Clare Boothe Luce Program.

\section{Smooth subpolytopes of $2\Delta_n$}
\label{sec:subpolytopes}

The purpose of this section is to study the polytopes $\Delta_n$,
$2\Delta_n$, $(2\Delta_n)_k$ and $\Delta_\ell\times\Delta_{n-\ell}$
appearing in the first four rows of {Table~\ref{table:dimdeg}}.  The
dimensions and degrees of the corresponding secant varieties will be
computed using determinantal methods.  We will also see that these
polytopes give essentially all smooth subpolytopes of $2\Delta_n$.

\subsection{Dimension and degree calculations} We begin with the
polytope $(2\Delta_n)_k$.  As in the introduction, this is defined to
be the convex hull of the lattice points of $2\Delta_n \setminus F$,
where $F$ is any $k$-dimensional face of $2\Delta_n$.  When $k = -1$,
$F$ is the empty face, so that $(2\Delta_n)_{-1} = 2\Delta_n$, and
when $k = n-1$, $F$ is a facet, so that $(2\Delta_n)_{n-1} =
\Delta_n$.  In the discussion that follows, we will usually exclude
these cases by requiring that $0 \le k \le n-2$.

The toric variety corresponding to $(2\Delta_n)_k$ is easy to describe.

\begin{prop}
\label{prop:dnrtoric}
Fix $k$ between $0$ and $n-2$, and let $Bl_k(\PP^n)$ denote the
blow-up of $\PP^n$ along a torus-invariant $k$-dimensional subspace.
Then, $Bl_k(\PP^n)$ is the toric variety $X_{(2\Delta_n)_k}$.
\end{prop}

\begin{proof}
Let $e_1,\dots,e_n$ be the standard basis of $\ZZ^n$.  The standard
fan for $\PP^n$ has cone generators $v_0,v_1,\dots,v_n$, where $v_1 =
e_1,\dots,v_n = e_n$ and $v_0 = -\sum_{i=1}^n e_i$.  Given $k$ between
$0$ and $n-2$, the star subdivison of the $(n-k)$-dimensional cone
$\sigma = \mathrm{Cone}(v_0,v_{k+2},\dots,v_n)$ is obtained by adding
the new cone generator
\[
v_{n+1} = v_0 + v_{k+2} + \dots + v_n = -e_1 - \dots - e_{k+1}
\]
and subdivding $\sigma$ accordingly.  As explained in
\cite[Prop.~1.26]{oda}, the toric variety of this new fan is the blow
up of $\PP^n$ along the $k$-dimensional orbit closure corresponding to
$\sigma$.  Thus we have $Bl_k(\PP^n)$.  The cone generators
$v_0,\dots, v_{n+1}$ give torus-invariant divisors $D_0,\dots,
D_{n+1}$ on $Bl_k(\PP^n)$.  One easily checks that $D = 2D_0 +
D_{n+1}$ is ample.

Using the description on p.\ 66 of \cite{fulton-t}, the polytope $P_D$
determined by $D$ is given by the $n+2$ facet inequalities
\begin{align*} 
x_1 + \dots + x_n &\le 2\\ 
x_i &\ge 0, \quad 1 \le i \le n\\ 
x_1 + \dots + x_{k+1} &\le 1.
\end{align*}
The inequalities on the first two lines define $2\Delta_n$ and the final inequality
removes the face corresponding to $\sigma$.  It follows that $P_D =
(2\Delta_n)_k$.
\end{proof}

We will see below that the ideal of $Bl_k(\PP^n)$ is determinantal.
We now turn our attention to some dimension and degree calculations.

\begin{thm}
\label{thm:dimdegspecial}
If $P$ is one of the polytopes
\[
\Delta_n,\ \ 2\Delta_n,\ \ (2\Delta_n)_k\ (0 \le k \le n-2),\ \
\Delta_\ell \times \Delta_{n-\ell}\ (1 \le \ell \le n-1),
\]
then the dimension and degree of $\Sec X_P$ are given by the first
four rows of {Table~\ref{table:dimdeg}}.
\end{thm}

\begin{proof} The case of $\Delta_n$ is trivial.  Turning our
attention to $2\Delta_n$, we let $r = \binom{n+2}{2}-1$.  The lattice
points in $2\Delta_n$ correspond to the $r+1$ monomials of degree 2 in
$x_0, \ldots, x_n$.  Order these monomials lexicographically so that
$\nu_2:\PP^n \to \PP^r$ is the map
\[
[x_0:\cdots:x_n] \mapsto [x_0^2: x_0x_1:\cdots:x_n^2].
\]  
If $z_0, \ldots, z_r$ are the homogeneous coordinates on $\PP^r$, then
the $2 \times 2$ minors of the $(n+1)\times(n+1)$ symmetric matrix
\[
M = \begin{pmatrix}
z_0 & z_1 & \cdots & z_n\\
z_1 & z_{n+1} & \cdots & z_{2n}\\
\vdots & \vdots && \vdots\\
z_n & z_{2n} & \cdots & z_r 
\end{pmatrix}
\]
vanish on $\nu_2(\PP^n)$.  As noted in Example~\ref{ex:veronese},
these minors generate the ideal of $\nu_2(\PP^n)$ and the $3 \times 3$
minors generate the ideal of $\Sec \nu_2(\PP^n)$.  The degree of this
determinantal ideal is also classical---see Example 14.4.14 of
\cite{fulton-i} for a reference.  This completes the proof for
$2\Delta_n$. 

Since $(2 \Delta_n)_k$ is constructed from $2 \Delta_n$ by removing a
face of dimension $k$, $(2 \Delta_n)_k$ has exactly $\binom{k+2}{2}$
fewer lattice points.  We will use the convention that these lattice
points correspond to the last $\binom{k+2}{2}$ monomials ordered
lexicographically.  The corresponding variables lie in the last $k+1$
rows and columns of $M$.

Let $M_{k+1}$ denote the matrix $M$ minus its last $k+1$ rows.  The
matrix $M_{k+1}$ is the \emph{partially symmetric} $(n-k) \times (n+1)$
matrix of Remark 2.5 (c) in \cite{conca}.  Since the $2 \times 2$
minors of $M_{k+1}$ vanish on $X_{(2 \Delta_n)_k}$, it follows that if
they generate a prime ideal defining a projective variety of dimension
$n$, then they generate the ideal of $X_{(2 \Delta_n)_k}$.  Moreover,
if the $2 \times 2$ minors of $M_{k+1}$ generate the ideal of $X_{(2
\Delta_n)_k}$, then the $3 \times 3$ minors of $M_{k+1}$ vanish on its
secant variety.  Thus, if the $3 \times 3$ minors generate a prime
ideal defining a projective variety of dimension $\dim \Sec X_{(2
\Delta_n)_k},$ then the
ideal they generate must be the ideal of the secant variety.  

From 
Example 3.8 of \cite{conca} we know that the dimension and degree of the projective variety cut out by the ideal generated by the $t \times t$ minors of $M_{k+1}$
  is what we desire for $t = 1, \ldots, n-k$, and from Remark 2.5 of \cite{conca} we know that the ideal is prime.  Moreover, $X_{(2\Delta_n)_{n-2}}$ is the rational normal scroll $S_{1, \ldots, 1,2},$ and by \cite{cj2}, we know that its
secant variety fills $\PP^{2n}.$  Therefore, we just need to show that the secant variety of $X_{(2\Delta_n)_k}$ also has dimension $2n$ for $0 \leq k < n-2.$  Note that $X_{(2\Delta_n)_k }$ is the image 
of $X_{(2\Delta_n)_j}$ via a linear projection for all $j < k.$  The variety $X_{(2\Delta_n)_{-1}}$ is $\nu_2(\PP^n).$    Therefore, applying Lemma 1.12 of \cite{cr} twice, tells us that $\dim \Sec X_{(2\Delta_n)_k } = 2n$ for all $k = 0, \ldots, n-2.$

Finally, $\Delta_\ell \times \Delta_{n-\ell}$ gives the Segre
embedding of $\PP^\ell \times \PP^{n-\ell}$.  In this case, the book \cite{fov} contains an explicit
description of the ideal of the secant variety in Example 1.3.6 (2). The 
variety $\Sec \PP^\ell \times \PP^{n-\ell}$ is
defined by the $3 \times 3$ minors of the generic $(\ell+1)\times
(n-\ell+1)$ matrix if $\ell+1, n-\ell+1 \ge 3$.  The degree of this
determinantal ideal was known to Giambelli and the formula appearing
in {Table~\ref{table:dimdeg}} is given in Example 19.10 of
\cite{harris}.  If either $\ell+1$ or $n+1-\ell$ is strictly less than
3, then $\PP^\ell \times \PP^{n-\ell}$ embeds into $\PP^{2n-1}$, and
its secant variety fills the ambient space (see Example 4.5.21 in \cite{fov}), so that the degree is $1$.
This completes the proof of the theorem.
\end{proof}

The determinantal ideals arising from the polytopes $2\Delta_n$, $(2
\Delta_n)_k$ and $\Delta_\ell \times \Delta_{n-\ell}$ in the above
proof have many beautiful combinatorial and computational properties.
Under the lexicographic term ordering (and other sufficiently nice
term orderings), the $t \times t$ minors contain a Gr\"obner basis for
the ideal they generate, and this ideal has a square-free initial
ideal.  Thus, Gr\"obner basis techniques can be used to show that the
minors generate radical ideals, and Stanley-Reisner techniques can be
used to compute degrees.  This can be proved using \cite{conca} and
the papers cited in \cite{conca}.

We also note that for these polytopes $P$, the ideal of $X_P$ is
generated by quadrics (the $2\times2$ minors) and the ideal of $\Sec
X_P$ is generated by cubics (the $3\times3$ minors).  Are there
other interesting polytopes with these properties?  Such questions
have been raised in \cite[Problem 5.15]{erss} in the context of
Jukes-Cantor binary models and in \cite[Conjecture 3.8]{vermeire},
which in the toric case asks whether the ideal of $\Sec X_{mP}$ is
generated by cubics for $m \gg 0$.

\subsection{Classification} The polytopes appearing in
Theorem~\ref{thm:dimdegspecial} fit naturally inside $2\Delta_n$.  We
will now prove that these are essentially all smooth polytopes with
this property.  

The proofs given here and in Section~\ref{sec:complete} require three
lemmas about smooth polytopes.  We begin with some defintions.

\begin{defin}
\label{def:edgelength}
If $\sigma$ is an edge of a lattice polytope $P$, then its \emph{edge
length} is $|\sigma\cap\ZZ^n|-1$, or alternatively, the normalized
length of $\sigma$.  We say that $P$ is a polytope of \emph{edge
length~1} if its edges all have length~1.
\end{defin}

In terms of the corresponding toric variety $X_P$, $P$ determines an
ample divisor $D_P$ and $\sigma$ determines a curve $C_\sigma$.  The
intersection product $D_P\cdot C_\sigma$ is the edge length of
$\sigma$.

\begin{defin}
\label{def:standardposition}
A vertex $v$ of a smooth polytope $P$ is in \emph{standard position}
if $v$ is the origin and its nearest lattice neighbors along the
edges are the standard basis $e_1,\dots,e_n$ of $\ZZ^n$.
\end{defin}

It is easy to see that any vertex of a smooth polytope can be moved
into standard position via an element of $AGL_n(\ZZ)$.  

We omit the elementary proof of our first lemma.

\begin{lem}
\label{lem:point11}
Let $P$ be a smooth polygon with a vertex at the origin in standard
position.  If one of edges at the origin has edge length $\ge 2$, then
$P$ contains the point $e_1+e_2$.
\end{lem}

\begin{lem}
\label{lem:dim3}
If $P$ is a smooth 3-dimensional polytope of edge length~1, then, up
to $AGL_3(\ZZ)$-equivalence, one of the following holds
\begin{align*}
&P = \Delta_3\\
&P = \Delta_1\times\Delta_2\\
&P\ \text{\rm contains}\ \Delta_1^3\\
&P\ \text{\rm contains}\ P_{1,2,2} =
\mathrm{Conv}(0,e_1,e_2,2e_1+e_3,2e_2+e_3,e_3). 
\end{align*}
\end{lem}

\begin{rmk} Even though $P_{1,2,2}$ has two edges of length~2,
it can be contained in a smooth 3-dimensional polytope of edge
length~1.
\end{rmk}

\begin{proof}
Put a vertex of $P$ in standard position and note that any facet of $P$
containing the origin is one of the following three types: (A) the
facet is a triangle $\Delta_2$; (B) it is a square $\Delta_1^2$; or
(C) it contains the lattice points
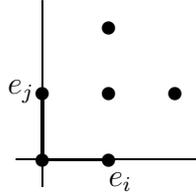
\[
\begin{picture}(100,70)
\put(15,10){\circle*{5}}
\put(15,35){\circle*{5}}
\put(40,60){\circle*{5}}
\put(40,10){\circle*{5}}
\put(65,35){\circle*{5}}
\put(40,35){\circle*{5}}
\put(5,10){\line(1,0){70}}
\put(15,0){\line(0,1){70}}
\put(40,0){$e_i$}
\put(2,35){$e_j$}
\thicklines
\put(15,10){\line(1,0){25}}
\put(15,10){\line(0,1){25}}
\end{picture}
\]
Now consider the three facets that meet at the origin.  One of
four things can happen: 
\begin{itemize}
\item Two facets are of type (A).  Then $P = \Delta_3$.
\item One facet is of type (A) and at least one is of type (B).  Then $P =
\Delta_1\times\Delta_2$.
\item Two facets are of type (C).  Then $P$ contains $P_{1,2,2}$.
\item Two facets are of type (B) and the other is not of type (A).  Then
$P$ contains $\Delta_1^3$.
\end{itemize}
The proofs of the first three bullets are elementary and are omitted.
For the fourth, we can assume that the two type (B) facets meet along
the edge connecting 0 to $e_1$.  The third facet meeting at $e_1$ is
of type (B) or (C), which easily implies that $P$ contains
$\Delta_1^3$.
\end{proof}

\begin{lem} 
\label{lem:3dimfaces}
Let $P$ be a smooth polytope of dimension $n \ge 3$ and edge length~1.
If every 3-dimensional face of $P$ is $AGL_3(\ZZ)$-equivalent to
$\Delta_3$ or $\Delta_1 \times \Delta_2$, then $P$ is
$AGL_n(\ZZ)$-equivalent to $\Delta_\ell\times\Delta_{n-\ell}$, where
$0 \le \ell \le n$.
\end{lem}

\begin{proof}
When we put $P$ into standard position, the edge length~1 hypothesis
implies that $0,e_1,\dots,e_n$ are vertices of $P$.  The vertex $e_1$
lies in $n$ facets of $P$, $n-1$ of which lie in coordinate
hyperplanes.  For the remaining facet $\Gamma$, its supporting
hyperplane can be defined by
\[
x_1 + a_2 x_2 + \dots + a_n x_n = 1, \quad a_i \in \QQ. 
\]
At $e_1$, the nearest neighbors along the edges consist of 0 and $n-1$
vertices lying in $\Gamma$.  Our hypothesis on $P$ implies that every
2-dimensional face of $P$ is either $\Delta_2$ or $\Delta_1^2$.  By
considering the 2-dimensional face contained in
$\mathrm{Cone}(e_1,e_i)$ for $i = 2,\dots,n$, we see that the $n-1$
nearest neighbors in $\Gamma$ are either $e_i$ or $e_i + e_1$.  We can
renumber so that the nearest neighbors in $\Gamma$ are
\[
e_2,\dots,e_\ell,e_{\ell+1}+e_1,\dots,e_n+e_1
\]
for some $1 \le \ell \le n$.  It follows that the supporting hyperplane
of $\Gamma$ is
\begin{equation}
\label{eq:1stsupp}
x_1+\dots+x_\ell = 1.
\end{equation}
If $\ell = n$, then $P = \Delta_n$ follows easily.  So now assume that
$\ell < n$ and pick $i$ between $2$ and $n-1$.  Consider the
3-dimensional face $F_i$ of $P$ lying in $\mathrm{Cone}(e_1, e_i,
e_n)$.  Depending on $i$, we get the following partial picture of
$F_i$:
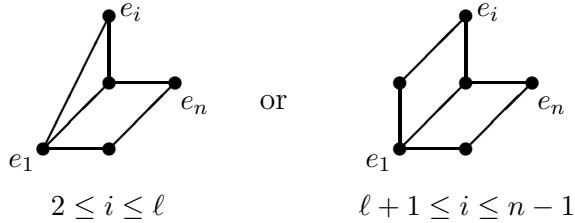
\[
\begin{picture}(235,85)(0,-15)
\put(15,10){\circle*{5}}
\put(40,35){\circle*{5}}
\put(40,10){\circle*{5}}
\put(65,35){\circle*{5}}
\put(40,60){\circle*{5}}
\put(97,25){or}
\put(18,-15){\small $2 \le i \le \ell$}
\put(134.5,-15){\small $\ell+1 \le i \le n-1$}
\put(150,10){\circle*{5}}
\put(150,35){\circle*{5}}
\put(175,60){\circle*{5}}
\put(175,10){\circle*{5}}
\put(175,35){\circle*{5}}
\put(200,35){\circle*{5}}
\put(67,25){\small $e_n$}
\put(44,60){\small $e_i$}
\put(2,3){\small $e_1$}
\put(202,25){\small $e_n$}
\put(179,60){\small $e_i$}
\put(137,3){\small $e_1$}
\thicklines
\put(15,10){\line(1,0){25}}
\put(15,10){\line(1,2){25}}
\put(40,35){\line(1,0){25}}
\put(40,10){\line(1,1){25}}
\put(15,10){\line(1,1){25}}
\put(40,35){\line(0,1){25}}
\put(150,10){\line(1,0){25}}
\put(150,10){\line(0,1){25}}
\put(150,10){\line(1,1){25}}
\put(150,35){\line(1,1){25}}
\put(175,10){\line(1,1){25}}
\put(175,35){\line(1,0){25}}
\put(175,35){\line(0,1){25}}
\end{picture}
\]
This picture shows that $F_i$ cannot be $\Delta_3$.  Hence by
assumption it is $AGL_3(\ZZ)$-equivalent to $\Delta_1 \times
\Delta_2$.  In terms of the picture, this means that the nearest
neighbors of the vertex $e_n$ of $F_i$, are 0, $e_1+e_n$, and
\begin{align*}
e_i+e_n\quad &\text{if}\ 2 \le i \le \ell\\
e_i\quad &\text{if}\ \ell+1 \le i \le n-1.
\end{align*}
Thus the nearest neighbors along the edges of $P$ at $e_n$ are
\[
0,e_1+e_n,\dots,e_\ell+e_n,e_{\ell+1},\dots,e_{n-1}.
\]
This gives the supporting hyperplane
\begin{equation}
\label{eq:2ndsupp}
x_{\ell+1}+\dots+x_n = 1.
\end{equation}
The combination of \eqref{eq:1stsupp} and \eqref{eq:2ndsupp} implies
easily that $P = \Delta_\ell \times \Delta_{n-\ell}$.
\end{proof}

We now prove a preliminary version of our classification of smooth
subpolytopes of $2\Delta_n$.

\begin{thm}
\label{thm:prelimclass}
Let $P$ be a smooth polytope of dimension $n$ in standard position at
the origin.  If $P \subset 2\Delta_n$, then $P$ is one of
\[
\Delta_n,\ \ 2\Delta_n,\ \ (2\Delta_n)_k\ (0 \le k \le n-2),\ \
\Delta_\ell \times \Delta_{n-\ell}\ (1 \le \ell \le n-1)
\]
up to $AGL_n(\ZZ)$-equivalence.
\end{thm}

\begin{rmk} In Section~\ref{sec:complete}, we will show that the
standard position hypothesis in
Theorem~\ref{thm:prelimclass} is unnecessary. 
\end{rmk}

\begin{proof}
The proof is trivial when $n = 1,2$.  So we will assume that $n \ge
3$.  For each $1 \le i \le n$, the edge starting from $0$ in direction
$e_i$ ends at either $e_i$ or $2e_i$.  Renumbering if necessary (which
can be done by an $AGL_n(\ZZ)$-equivalence), we can assume that
\begin{equation}
\label{eq:somevertP}
e_1,\dots,e_{k+1},2e_{k+2},\dots,2e_n
\end{equation}
are vertices of $P$.  The case $k = -1$ corresponds to $P =
2\Delta_n$.  Now suppose that $0 \le k \le n-2$, so that $2e_n$ is a
vertex of $P$.  Given $1 \le i < j \le k+1$, we claim that $e_i+e_j
\notin P$.  To prove this, assume $e_i+e_j \in P$ and consider the
3-dimensional face $F$ of $P$ determined by the vertices
$0,e_i,e_j,2e_n$:
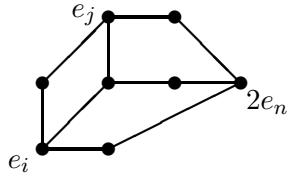
\[
\begin{picture}(100,68)(0,-3)
\put(15,10){\circle*{5}}
\put(15,35){\circle*{5}}
\put(40,35){\circle*{5}}
\put(40,10){\circle*{5}}
\put(65,35){\circle*{5}}
\put(90,35){\circle*{5}}
\put(40,60){\circle*{5}}
\put(65,60){\circle*{5}}
\put(92,25){\small $2e_n$}
\put(26,60){\small $e_j$}
\put(2,3){\small $e_i$}
\thicklines
\put(15,10){\line(1,0){25}}
\put(65,60){\line(1,-1){25}}
\put(15,10){\line(0,1){25}}
\put(15,35){\line(1,1){25}}
\put(40,35){\line(1,0){50}}
\put(40,10){\line(2,1){50}}
\put(15,10){\line(1,1){25}}
\put(40,35){\line(0,1){25}}
\put(40,60){\line(1,0){25}}
\end{picture}
\]
Note that $e_i+e_n, e_j + e_n \in P$ by Lemma~\ref{lem:point11}.
Since $P \subset 2\Delta_n$, we have $F \subset 2\Delta_3$, so that
the above picture shows all lattice points of $F$.  The convex hull of
these points is not smooth, giving the desired contradiction.  It
follows that
\[
P \subset 2\Delta_n \setminus \mathrm{Conv}(2e_1,\dots,2e_{k+1}),
\]
which easily implies that $P \subset (2\Delta_n)_k$.  For the opposite
inclusion, note that since \eqref{eq:somevertP} consists of lattice
points of $P$, Lemma~\ref{lem:point11} implies that $P$ also contains
the lattice point $e_i+e_j$ whenever $1 \le i \le k+1$ and $k+2 \le j
\le n$.  Since $(2\Delta_n)_k$ is the convex hull of these points
together with \eqref{eq:somevertP}, we obtain $(2\Delta_n)_k \subset
P$.

It remains to consider the case when $P$ contains none of
$2e_1,\dots,2e_n$.  This implies that $P$ has edge length~1.  Let $F$
be a 3-dimensional face of $P$.  Of the four possibilities for $F$
listed in Lemma~\ref{lem:dim3}, note that $F$ cannot contain a
configuration $AGL_3(\ZZ)$-equivalent to $P_{1,2,2}$, since
$P_{1,2,2}$ contains two parallel segments with $3$ lattice points,
which is impossible in $2\Delta_n$.  Similarly, $2\Delta_n$ cannot
contain a configuration $AGL_3(\ZZ)$-equivalent to $\Delta_1^3$.  This
is because such a ``cube'' configuration would give an affine relation
\[
(v_1 - v_0) + (v_2-v_0) + (v_3-v_0) = v_4-v_0,
\]
where $v_0$ and $v_4$ are opposite vertices of the cube and
$v_1,v_2,v_3$ are vertices of the cube nearest to $v_0$.  This gives
the relation
\[
v_1+v_2+v_3 = 2v_0 + v_4,
\]
which cannot occur among distinct lattice points of $2\Delta_n$ (we
omit the elementary argument).  It follows from Lemma~\ref{lem:dim3}
that $F$ is $AGL_3(\ZZ)$-equivalent to $\Delta_3$ or
$\Delta_1\times\Delta_2$.  Since $F$ is an arbitary 3-dimensional face
of $P$, Lemma~\ref{lem:3dimfaces} implies that $P$ is
$AGL_n(\ZZ)$-equivalent to $\Delta_\ell \times \Delta_{n-\ell}$ for
some $0 \le \ell \le n$.  The proof of the theorem is now complete.
\end{proof}

\section{The degree formula for the secant variety}
\label{sec:degformula}

One of the key ingredients in the proof of Theorem~\ref{thm:dimension}
is a formula relating the degree of $\Sec X$ to the Chern classes of
the tangent bundle $T_X$ and the line bundle $\mathcal{O}_X(1)$.  We
review this formula and then interpret it for toric surfaces and
3-folds.

\subsection{The abstract join and the secant variety}
\label{sec:absjoin}
Given $X \subset \PP^r$ of dimension $n$, consider $\PP^{2r+1}$ with
coordinates $x_0, \dots, x_r, y_0, \dots, y_r$.  The \emph{abstract
join} $J(X,X)$ of $X$ with itself is the set of all points in
$\PP^{2r+1}$ of the form $[\lambda x: \mu y]$ with $[x], [y] \in X$
and $[\lambda:\mu] \in \PP^1$.  Then $J(X,X)$ has dimension $2n+1$ and
degree $(\deg X)^2$ (see \cite{fov}).

The linear projection 
\[
\phi:\PP^{2r+1} \dashrightarrow \PP^r
\]
given by $\phi([x:y]) = [x-y]$ is defined away from the subspace of
$\PP^{2r+1}$ defined by the vanishing of $x_i-y_i$ for $i =
0,\dots,r$.  This induces a rational map $J(X,X) \dashrightarrow \Sec
X$ with base locus the diagonal embedding of $X$ in the abstract join.
We have the following well-known result:

\begin{thm}[Theorem 8.2.8 in \cite{fov}]
\label{thm:degreemap}
If $X$ is smooth, then 
\[
\deg \Sec X \deg \phi = (\deg X)^2 - \sum_{i = 0}^{n}
\binom{2n+1}{i} \int_X c(T_X)^{-1} \cap c_1(\mathcal{L})^i, 
\]
where $\mathcal{L} = \mathcal{O}_X(1)$ and $\deg \phi = 0$ when $\dim
\Sec X < 2n+1$.
\end{thm}

This is sometimes called the \emph{double point formula} since a
general linear projection $X \to \PP^{2n}$ has $\tfrac12\deg \Sec
X \deg \phi$ double points (see Corollary 8.2.6 of \cite{fov}, for
example).  For surfaces, the double point formula was first discovered
by Severi.  See \cite{fov,fulton-i} for more general double point
formulas and further references.

The following proposition explains the geometric meaning of $\deg
\phi$.  

\begin{prop}
\label{prop:degphi}
If $X \subset \PP^r$ is a variety of dimension $n > 1,$ then:
\begin{enumerate}
\item $\dim \Sec X = 2n+1$ if and only if a general point of $\Sec X$
lies on at most finitely many secant lines of $X$.
\item $\dim \Sec X = 2n+1$ implies that $\deg
\phi$ is $2$ times the number of secant lines of $X$ through a general
point of $\Sec X$.
\end{enumerate}
\end{prop}

\begin{proof}
Both of the claims are probably well-known to experts, but we could not find  
explicit proofs in the literature. Two different ways of constructing an ``abstract'' secant variety
which maps to $\Sec X$ may be found in \cite{zak} and \cite{cr} and 
one may treat the results here from those points of view as well.
We include a proof for completeness using the setup of \cite{fov}. 

If $\dim \Sec X = 2n+1$, then the rational map $\phi : J(X,X)
\dashrightarrow \Sec X$ is generically finite.  This implies
that a general point of $\Sec X$ lies on at most finitely many secant
lines of $X$.  Conversely, suppose that this condition is satisfied
and $\dim \Sec X < 2n+1$.  Then $\phi^{-1}(z)$ is infinite for $z \in
\Sec X$, so that we can find infinitely many distinct pairs $p_i \ne
q_i$ in $X$ such that $z$ lies on $\overline{p_iq_i}$.  Since there
are only finitely many such lines through $z$, one of them must
contain infinitely many points of $X$ and hence lies in $X$.  This
line also contains $z$, so that $z \in X$, and then $X = \Sec X$
follows.  One easily concludes that $X$ is a linear subspace.  Since
$X$ has dimension $> 1$, a general point of $\Sec X = X$ lies on
infinitely many secant lines, a contradiction.

For the second assertion, first observe that the trisecant variety of
$X$ (the closure of the union of secant lines meeting $X$ in $\ge 3$
points) has dimension $\le 2n$ by Corollary 4.6.17 of \cite{fov}.
Since $\dim \Sec X = 2n+1$, it follows that at a general point $z \in
\Sec X$, any secant line of $X$ through $z$ meets $X$ in exactly two
points, say $p \ne q$.  Writing $z = \lambda p + \mu q$ gives distinct
points $[\lambda p:-\mu q] \ne [\mu q:-\lambda p]$ in $\phi^{-1}(z)$.
Thus the cardinality of $\phi^{-1}(z)$ is twice the number of secant
lines through $z$, as claimed. 
\end{proof}

Here is an immediate corollary of Proposition~\ref{prop:degphi}.

\begin{cor}
\label{cor:degphi}
If $X \subset \PP^r$ has dimension $n > 1$, then $\deg \phi = 2$ if
and only if there is a unique secant line of $X$ through a general
point of $\Sec X$.  Furthermore, $\dim \Sec X = 2n+1$ when either of
these conditions are satisfied.
\end{cor}

Theorem~\ref{thm:degreemap} gives a formula for $\deg \Sec X$ when
$\deg \phi = 2$.  This explains the $\frac12$ appearing in
Theorem~\ref{thm:dimension} and Corollaries~\ref{thm:sfdegree}
and~\ref{thm:3folddegree}.

Here is an example of Theorem \ref{thm:degreemap} and
Corollary~\ref{cor:degphi} that will be useful in the proof of Theorem
\ref{thm:dimension}.

\begin{ex}
\label{ex:rns}
Consider the rational normal scroll $X = S_{d_1,\dots,d_n}$, where
$d_i \ge 1$ and $d = \sum_{i=1}^n d_i$.  Since $X$ is a projective
bundle over $\PP^1$, its Chow ring is well-known, making it easy to
compute the formula in Theorem \ref{thm:degreemap}.  This computation
appears in \cite{ohno}, with the result
\[
\deg \Sec X \deg \phi = d^2 - (2n+1)d + n(n+1).
\]
Catalano-Johnson proved in \cite{cj} that a general point of $\Sec X$
lies on a unique secant line of $X$ when $d \ge n+2$.  Thus $\deg \phi
= 2$, so that
\[
\deg \Sec X = \tfrac12\big(d^2 - (2n+1)d + n(n+1)\big), \quad d \ge
n+2.
\]
If $d = n, n+1$, this formula gives zero, so that $\dim \Sec X < 2n+1$
by Theorem~\ref{thm:degreemap}.  The cases $d = n, n+1$ correspond to
the polytopes $\Delta_1 \times \Delta_{n-1}$, and $(2\Delta_n)_{n-2}$,
respectively.  Since the secant variety of a rational normal scroll
always has the expected dimension, this explains why these polytopes
don't appear in the statement of Corollary~\ref{cor:dimension}.
\end{ex}

\subsection{Dimensions 2 and 3}
\label{sec:dim2and3}
In low dimensions, the degree formula of Theorem \ref{thm:degreemap}
can be expressed quite succinctly.  The purpose of this section is to
prove the following two theorems.

\begin{thm}\label{thm:sfdegreegen}
If $P$ is a smooth lattice polygon in $\RR^2$, then
\[
\deg \Sec X_P \deg \phi = d^2-10d+5B+2V-12,
\]
where $d$ is the area of $P,$ $B$ is the number of lattice points on
the boundary of $P$, and $V$ is the number of vertices of $P$.
\end{thm}

\begin{thm}\label{thm:3folddegreegen}
If $P$ is a smooth 3-dimensional lattice polytope in $\RR^3$, then
\[
\deg \Sec X_P \deg \phi = d^2-21d+c_1^3+8V+14E-84I-132,
\]
where $d$ is the volume of $P$, $E$ is the number of lattice points on
the edges of $P$, $V$ is the number of vertices of $P$, and $I$ is the
number of interior lattice points.  Also, $c_1 = c_1(T_{X_P})$.
\end{thm}

Both results follow from Theorem \ref{thm:degreemap} via applications
of Riemann-Roch and the theory of Ehrhart polynomials.  We begin with
some useful facts about the total Chern class $c(T_X)$ and Todd class
$\mathrm{td}(X)$ from pp.\ 109--112 in \cite{fulton-t}.  Let $c_i =
c_i(T_X)$ and $H = c_1(\mathcal{L})$, where $\mathcal{L}$ is the line
bundle coming from $P$.  Then:
\begin{equation}
\label{eq:todd}
\begin{aligned}
\mathrm{td}(X) &= 1+\tfrac12 c_1 + \tfrac{1}{12}(c_1^2+c_2) +
\tfrac{1}{24}c_1c_2+ \cdots\\   
c(T_X)^{-1} &= 1-c_1+(c_1^2-c_2) + (2c_1c_2-c_1^3-c_3)+ \cdots.
\end{aligned}
\end{equation}
Furthermore, if $X =X_P$ is the smooth toric variety of the polytope
$P$ and $V(F)$ is the orbit closure corresponding to the face $F$ of
$P$, then:
\begin{equation}
\label{eq:chern}
\begin{aligned}
c(T_X) &= \prod_{\dim F = n-1} (1+[V(F)])\\
c_i &= \sum_{\dim F = n-i} [V(F)]\\
H^i \, \cap [V(F)] &= \mathrm{Vol}_{i}(F),
\end{aligned}
\end{equation}
where $\mathrm{Vol}_{i}(F)$ is the normalized volume of the
$i$-dimensional face $F$.

\begin{proof}[Proof of Theorem \ref{thm:sfdegreegen}]
Since $\dim X =2$ and $\deg X = H^2 = d$, the right-hand side of
the formula in Theorem \ref{thm:degreemap} becomes
\[
d^2 - (c_1^2-c_2)+5Hc_1-10d.
\]
Since $c_1^2+c_2 = 12$ by Noether's formula, this simplifies to
\[
d^2 - 12+2c_2+5Hc_1-10d.
\]
By \eqref{eq:chern}, $c_2$ is the sum of the torus-fixed points
corresponding to the vertices of $P$.  Thus $c_2 = V$, the number of
vertices of $P$.  Furthermore, \eqref{eq:chern} also implies that $Hc_1$
is the perimeter of $P$, which is the number $B$ of lattice points on
the boundary of $P$.  The desired formula follows.
\end{proof}

\begin{proof}[Proof of Theorem \ref{thm:3folddegreegen}]
The formula of Theorem \ref{thm:degreemap} reduces to
\[
d^2 - (2c_1c_2-c_1^3-c_3) -7 H (c_1^2-c_2) +21 H^2c_1-35 d,
\]
where $d$ again denotes $\deg X$.

Let $S$ denote the surface area of $P$, $\mathcal{E}$ denote the
perimeter (i.e., the sum of the lengths of the edges of $P$), and $V$
denote the number of vertices of $P$.  Using \eqref{eq:chern} as
before, one sees that $c_3= V$, $Hc_2 = \mathcal{E}$, and $H^2c_1 =
S$.  Also, using \eqref{eq:todd} together with the fact that
$\mathrm{td}_3(X) = [x]$ for any $x \in X$, we obtain $c_1c_2 = 24$.
Thus, we have
\[
d^2 - 48 + c_1^3 + V - 7 Hc_1^2 + 7 \mathcal{E} + 21S - 35 d.
\]
However, the Todd class formula from \eqref{eq:todd} and Riemann-Roch
tell us that the number of lattice points in $P$ is
\[
\ell(P)=1 + \tfrac1{12}(Hc_1^2+\mathcal{E}) + \tfrac14S + \tfrac16 d,
\]
as explained in \cite[Sec.\ 5.3]{fulton-t}.  Solving for $Hc_1^2$ and
substituting into the above expression, we obtain
\[
d^2-21d+c_1^3 + \underbrace{V + 14\mathcal{E}}_{\alpha} +
\underbrace{42S-84\ell(P)+36}_{\beta}. 
\]
It remains to understand the quantities $\alpha$ and $\beta$. 

Letting $E$ denote the number of lattice points on the edges of $P$,
we can rewrite $\alpha$ as $8V + 14E$ since every vertex of $P$ lies
on exactly three edges of $P$ by smoothness.

Next let $B$ (resp.\ $I$) denote the number of boundary (resp.\
interior) lattice points of $P$, so that $\ell(P) = B + I$.  By
Ehrhart duality, we have
\[
I = (-1)^3\big(1 - \tfrac1{12}(Hc_1^2+\mathcal{E}) + \tfrac14S -
\tfrac16 d\big),
\]
which allows us to rewrite $\beta$ as $-84I - 132$.  The desired
expression follows immediately.
\end{proof}

\section{Counting secant lines}
\label{sec:complete}

To complete the proofs of our main results, we need to study the
secant variety of a toric variety coming from the last row of
{Table~\ref{table:dimdeg}}.

\subsection{A uniqueness theorem} \label{sec:key} 
Here is our result.
 
\begin{thm} 
\label{thm:xp}
If a smooth $n$-dimensional polytope $P$ is not
$AGL_n(\ZZ)$-equivalent to any of the polytopes
\[
\Delta_n,\ \ 2\Delta_n,\ \ (2\Delta_n)_k\ (0 \le k \le n-2),\ \
\Delta_\ell \times \Delta_{n-\ell}\ (1 \le \ell \le n-1),
\]
then a unique secant line of $X_P$ goes through a general point of
$\Sec X_P$.
\end{thm}

Before beginning the proof of Theorem~\ref{thm:xp}, we need some
lemmas.

\begin{lem}
\label{lem:projection}
Let $X \subset \PP^r$ be a variety of dimension $n > 1$ and fix a
linear projection $\PP^r \dashrightarrow \PP^s$ such that $X$ is not
contained in the center of the projection.  Let $Y \subset \PP^s$ be
the closure of the image of $X$, so that we have a projection $\pi : X
\dashrightarrow Y$.  If $\pi$ is birational and a general point of
$\Sec Y$ lies on a unique secant line of $Y$, then $\dim \Sec X =
2n+1$ and a general point of $\Sec X$ lies on a unique secant line of
$X$.
\end{lem}

\begin{proof}
By Proposition~\ref{prop:degphi}, our hypothesis on $Y$ implies that
$\Sec Y$ has dimension $2n+1$, and then $\dim \Sec X = 2n+1$ since
$\pi : X \dashrightarrow Y$ induces a dominating map $\Sec X
\dashrightarrow \Sec Y$ .

Now suppose that a general point $z \in \Sec X$ lies on the secant
lines $\overline{pq}$ and $\overline{p'q'}$ of $X$.  These map to
secant lines of $Y$ through $\pi(z)$, which coincide by hypothesis.
Arguing as in the proof of Proposition~\ref{prop:degphi}, we can
assume that this secant line meets $Y$ at exactly two points.
Switching $p'$ and $q'$ if necessary, we get $\pi(p) = \pi(p')$ and
$\pi(q) = \pi(q')$.  Since $\pi$ is generically 1-to-1 on $X$, we
conclude that $\overline{pq} = \overline{p'q'}$.  Thus a general point
of $\Sec X$ lies on a unique secant line of $X$.
\end{proof}

We will prove Theorem~\ref{thm:xp} by applying
Lemma~\ref{lem:projection} to projections constructed from carefully
chosen subsets $A = \{u_0,\dots,u_s\} \subset P\cap\ZZ^n$.  The
characters $\chi^{u_i}$ give a rational map $X \dashrightarrow \PP^s$
defined by $x \mapsto [\chi^{u_0}(x):\cdots:\chi^{u_s}(x)]$, where $X$
is the abstract toric variety of the normal fan of $P$.  The closure
of the image in $\PP^s$ is denoted $X_A$.  Note that $X_A = X_P$ when
$A = P \cap \ZZ^n$.

We can view $X_A$ as a projection of $X_P$ as follows. By definition,
$X_P$ is the closure of the image of the rational map $X
\dashrightarrow \PP^r$, $r = |P\cap\ZZ^n| - 1$.  Labeling the
coordinates of $\PP^r$ using the lattice points of $P\cap\ZZ^n$, we
obtain a projection
\begin{equation}
\label{eq:projectionpa}
\pi:\PP^r \dashrightarrow \PP^s
\end{equation}
by projecting onto the linear subspace defined by the coordinates
corresponding to $A \subset P\cap\ZZ^n$.  This induces the projection
$\pi:X_P \dashrightarrow X_A$.

Some of the subsets $A$ that we will use come from the toric interpretation
of the rational normal scrolls from the introduction.  Let
$d_1,\dots,d_n$ be positive integers and label the vertices of the
unit simplex $\Delta_{n-1} \subset \RR^{n-1}$ as $v_1 =
e_1,\dots,v_{n-1} = e_{n-1}, v_n = 0$.  Define $A_{d_1,\dots,d_n}
\subset \ZZ^{n-1} \times \ZZ$ by
\[
A_{d_1,\dots,d_n} = \bigcup_{i=1}^n \big\{v_i + a_i e_n \mid a_i \in
\ZZ, 0 \le a_i \le d_i\big\},
\]
and let $P_{d_1,\dots,d_n}$ be the convex hull of $A_{d_1,\dots,d_n}$
in $\RR^{n-1}\times\RR$.  It is straighforward to verify that
$P_{d_1,\dots,d_n}$ is a smooth polytope with lattice points
\[
A_{d_1,\dots,d_n} = P_{d_1,\dots,d_n} \cap (\ZZ^{n-1}
\times \ZZ).
\]
One can easily show that the toric variety $X_{P_{d_1,\dots,d_n}}$ is
the rational normal scroll $S_{d_1,\dots,d_n}$ defined in
Example~\ref{ex:scroll}.  The result of \cite{cj} mentioned in
Example~\ref{ex:rns} implies the following lemma.

\begin{lem}
\label{lem:goodAs}
If $A$ is one of the two sets
\begin{align*}
&A_{1,\dots,1,1,3} \subset \ZZ^n\quad \text{\rm (there are $n-1$
    1's, $ n \ge 1$)}\\ 
&A_{1,\dots,1,2,2} \subset \ZZ^n\quad \text{\rm (there are $n-2$
    1's, $ n \ge 2$)},
\end{align*}
then a general point of $\Sec X_A$ lies on a unique secant line of $X_A$.
\end{lem}

For polytopes of edge length 1 and dimension $> 3$, we will use the
following configurations of lattice points.

\begin{lem}
\label{lem:moregoodA}
Let $B \subset \ZZ^3\times\{0\} \subset \ZZ^3 \times \ZZ^{n-3}$ be a
set of 8 lattice points such that $B$ affinely generates
$\ZZ^3\times\{0\}$ and a general point of $\Sec X_B$ lies on a unique
secant line of $X_B$.  For each $i = 4,\dots,n$, pick $u_i \in
\{e_1,e_2\}$ and define
\[
A = B \cup \{e_i, e_i+u_i \mid i = 4,\dots,n\} \subset \ZZ^3
\times \ZZ^{n-3}.
\]
Note that $A$ has $2n+2$ points.  If $X_A \subset \PP^{2n+1}$ is the
corresponding toric variety, then a general point of $\Sec X_A =
\PP^{2n+1}$ lies on a unique secant line of $X_A$.
\end{lem}

\begin{proof}
Let $t_1,\dots,t_n$ be torus variables corresponding to
$e_1,\dots,e_n$, and let $m_1,\dots,m_8$ be the lattice points of $B$.
A point $\mathbf{t} = (t_1,\dots,t_n)$ on the torus $(\CS)^n$ maps to
\[
\psi(\mathbf{t}) = (\mathbf{t}^{m_1},\dots,\mathbf{t}^{m_8},t_4,t_4
\mathbf{t}^{u_4},\dots,t_n,t_n\mathbf{t}^{u_n}) \in \PP^{2n+1}.
\]
We claim that the map $(\CS)^n \times (\CS)^n \times \CS
\dashrightarrow \Sec X_A$ defined by
\[
(\mathbf{t},\mathbf{s},\gamma) \mapsto \gamma \psi(\mathbf{t}) +
(1-\gamma) \psi(\mathbf{s})
\]
is generically 2-to-1.  To prove this, suppose that 
\begin{equation}
\label{eq:samept}
\gamma \psi(\mathbf{t}) + (1-\gamma) \psi(\mathbf{s}) = \gamma'
\psi(\mathbf{t}') + (1-\gamma') \psi(\mathbf{s}').
\end{equation}
We need to show that for $\mathbf{t}, \mathbf{s}, \gamma$ generic,
\eqref{eq:samept} implies that
\[
(\mathbf{t}',\mathbf{s}',\gamma') = (\mathbf{t},\mathbf{s},\gamma) \
\text{or} \ (\mathbf{s},\mathbf{t},1-\gamma).
\]
By projecting onto the first 8 coordinates and using our hypothesis on
$B$ and being careful to avoid trisecants and tangents, we obtain
\begin{equation}
\label{eq:first3}
(t_1',t_2',t_3',s_1',s_2',s_3',\gamma') =
(t_1,t_2,t_3,s_1,s_2,s_3,\gamma)
\end{equation}
by switching $\mathbf{t}$ and $\mathbf{s}$ if necessary.

Now fix $i$ between 4 and $n$.  Since $u_i \in \{e_1,e_2\}$,
\eqref{eq:first3} implies that $\mathbf{t}^{u_i} =
\mathbf{t}'{}^{u_i}$ and $\mathbf{s}^{u_i} = \mathbf{s}'{}^{u_i}$.
Since we also know that $\gamma = \gamma'$, comparing the coordinates
of \eqref{eq:samept} corresponding to $e_i, e_i + u_i \in B$ gives the
equations
\begin{align*}
\gamma t_i + (1-\gamma)s_i &= \gamma t_i' + (1-\gamma)s_i' \\
\gamma t_i \mathbf{t}^{u_i} + (1-\gamma)s_i \mathbf{s}^{u_i} &= \gamma
t_i'\mathbf{t}^{u_i} + (1-\gamma)s_i' \mathbf{s}^{u_i},
\end{align*}
which can be rewritten as
\begin{align*}
\gamma (t_i - t_i')+ (1-\gamma)(s_i - s_i') &= 0\\ \gamma
\mathbf{t}^{u_i} (t_i - t_i')+ (1-\gamma) \mathbf{s}^{u_i} (s_i -
s_i') &= 0.
\end{align*}
The coefficent matrix of this $2\times2$ system of homogeneous
equations has determinant
\[
\gamma(1-\gamma)(\mathbf{s}^{u_i} - \mathbf{t}^{u_i}),
\]
which is nonzero for generic $(\mathbf{t},\mathbf{s},\gamma)$.  Thus
$t_i' = t_i$ and $s_i' = s_i$ for $i = 4,\dots,n$.  It follows that
that the map $(\CS)^n \times (\CS)^n \times \CS \dashrightarrow \Sec
X_A$ defined above is generically 2-to-1.
\end{proof}

We are now ready to prove Theorem~\ref{thm:xp}.

\begin{proof}[Proof of Theorem~\ref{thm:xp}] 
The case of dimension 1 is trivial since $P = d\Delta_1$ gives a
rational normal curve of degree $d$ embedded in $\PP^d$.  It is
well-known (see \cite{cj} or Proposition~8.2.12 of \cite{fov}) that a
unique secant line passes through a general point of $\Sec X_P$ when
$d \ge 3$.

Now assume that $\dim P \ge 2$.  We will consider five cases,
depending on the maximum edge length and dimension of $P$.

\medskip

\noindent{\bf Case 1.} $P$ has an edge $E$ of length $\ge 3$.  Put a
vertex of $E$ in standard position so that $E$ contains $e_n$.  Now
fix $i$ between $1$ and $n-1$ and consider the 2-dimensional face
$F_i$ of $P$ lying in $\mathrm{Cone}(e_i,e_n)$.  Then $F_i$ is a
smooth polygon in standard position containing the points
\[
0, e_i, e_n, 2e_n, 3e_n.
\]
By Lemma~\ref{lem:point11}, we conclude that $e_i + e_n \in F$.  Hence
$P$ contains the $2n+2$ points
\[
A = \{0, e_i, e_i+ e_n, e_n, 2e_n, 3e_n \mid i = 1,\dots,n-1\}.
\]
This is the set $A = A_{1,\dots,1,1,3}$ defined earlier.  By
Lemma~\ref{lem:goodAs}, a general point of $\Sec X_A$ lies in a unique
secant line of $X_A$.

The embeddings of $X_P$ and $X_A$ given by $P\cap\ZZ^n$ and $A$
respectively share the points $0,e_1,\dots,e_n$.  It follows easily
that the linear projection $X_P \dashrightarrow X_A$ given by
projecting onto the coordinates corresponding to points of $A$ is the
identity on the tori of $X_P$ and $X_A$ and hence is birational.  By
Lemma~\ref{lem:projection}, we conclude that a general point of $\Sec
X_P$ lies in a unique secant line of $X_P$.

\medskip

\noindent{\bf Case 2.} $P$ has an edge $E$ of length 2 but no edges of
length $\ge 3$.  As in Case 1, put a vertex of $E$ in standard
position so that $E$ contains $e_n$.  Now fix $i$ between $1$ and
$n-1$ and consider the 2-dimensional face $F_i$ of $P$ lying in
$\mathrm{Cone}(e_i,e_n)$.  Using Lemma~\ref{lem:point11} as in Case 1
shows that $F_i$ contains the points
\[
0, e_i, e_i + e_n, e_n, 2e_n.
\]
If for some $i$ the face $F_i$ also contains $e_i + 2e_n$, then we can
relabel so that $i = n-1$.  Hence $P$ contains the points
\[
A = \{0, e_i, e_i+ e_n, e_{n-1}, e_{n-1}+ e_n, e_{n-1}+2e_n, e_n,
2e_n, \mid i = 1,\dots,n-2\}.
\]
This is the set $A = A_{1,\dots,1,2,2}$ defined earlier in the
section.  Then using Lemmas~\ref{lem:projection} and~\ref{lem:goodAs}
as in Case 1 implies that a general point of $\Sec
X_P$ lies in a unique secant line of $X_P$.

It remains to show that $e_i + 2e_n \notin F_i$ for all $i =
1,\dots,n-1$ cannot occur.  If this were to happen, it is easy to see
that $F_i$ is either
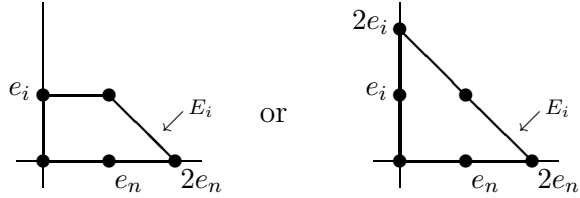
\[
\begin{picture}(235,70)
\put(15,10){\circle*{5}}
\put(15,35){\circle*{5}}
\put(40,10){\circle*{5}}
\put(40,35){\circle*{5}}
\put(65,10){\circle*{5}}
\put(5,10){\line(1,0){70}}
\put(15,0){\line(0,1){70}}
\put(97,25){or}
\put(150,10){\circle*{5}}
\put(150,35){\circle*{5}}
\put(150,60){\circle*{5}}
\put(175,10){\circle*{5}}
\put(175,35){\circle*{5}}
\put(200,10){\circle*{5}}
\put(140,10){\line(1,0){70}}
\put(150,0){\line(0,1){70}}
\put(67,0){\small $2e_n$}
\put(42,0){\small $e_n$}
\put(2,35){\small $e_i$}
\put(202,0){\small $2e_n$}
\put(177,0){\small $e_n$}
\put(137,35){\small $e_i$}
\put(132,60){\small $2e_i$}
\put(60,23){\tiny $\swarrow$}
\put(70,28){\tiny $E_i$}
\put(195,23){\tiny $\swarrow$}
\put(205,28){\tiny $E_i$}
\thicklines
\put(15,10){\line(1,0){50}}
\put(15,35){\line(1,0){25}}
\put(15,10){\line(0,1){25}}
\put(40,35){\line(1,-1){25}}
\put(150,10){\line(1,0){50}}
\put(150,10){\line(0,1){50}}
\put(150,60){\line(1,-1){50}}
\end{picture}
\]
Since $P$ is smooth, $n$ facets of $P$ meet at the vertex $2e_n$,
$n-1$ of which lie in coordinate hyperplanes.  The remaining facet has
a normal vector $\nu$ perpendicular to the edge $E_i$ of $F_i$
indicated in the above picture.  Thus $\nu\cdot(-e_i+e_n) = 0$ for $i
= 1,\dots,n-1$, which easily implies that the supporting hyperplane of
this facet is defined by $x_1 + \dots + x_n = 2$.  This hyperplane and
the coordinate hyperplanes bound $2\Delta_n$, and it follows that $P
\subset 2\Delta_n$.  Since $P$ is in standard position,
Theorem~\ref{thm:prelimclass} implies that $P$ is
$AGL_n(\ZZ)$-equivalent to one of the polytopes listed in the
statement of the theorem, a contradiction.

\medskip

\noindent{\bf Case 3.} $P$ has only edges of length 1 and dimension 2.
Put a vertex of $P$ in standard position and note that $e_1$ and $e_2$
are also vertices of $P$.  Since $P$ is smooth, the vertex $e_1$ of
$P$ must lie on an edge containing $a e_1 + e_2$, $a \ge 0$.  If $a =
0$ or 1, then one easily sees that $P$ is contained in $2\Delta_2$.
On the other hand, if $a \ge 3$, then $P$ contains
\[
A = \{0, e_1, e_2, e_1+e_2, 2e_1+e_2, 3e_1+e_2\},
\]
which equals $A_{1,3}$ up to $AGL_2(\ZZ)$-equivalence.  Using
Lemmas~\ref{lem:projection} and~\ref{lem:goodAs} as usual, we conclude
that a general point of $\Sec X_P$ lies in a unique secant line of $X_P$.  Finally, if
$a=2$, then $P$ is not contained in $2\Delta_2$.  Applying a similar
analysis to the edges emanating from the vertex $e_2$, one sees that
$P$ either contains $A_{1,3}$ up to $AGL_2(\ZZ)$-equivalence or has the
four solid edges pictured as follows:
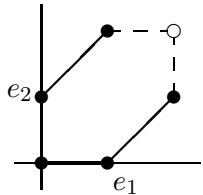
\[
\begin{picture}(100,70)
\put(15,10){\circle*{5}}
\put(15,35){\circle*{5}}
\put(40,60){\circle*{5}}
\put(40,10){\circle*{5}}
\put(65,35){\circle*{5}}
\put(5,10){\line(1,0){70}}
\put(15,0){\line(0,1){70}}
\put(42,0){\small $e_1$}
\put(2,35){\small $e_2$}
\put(40,60){\line(1,0){5}}
\put(50,60){\line(1,0){5}}
\put(60,60){\line(1,0){2.5}}
\put(65,35){\line(0,1){5}}
\put(65,45){\line(0,1){5}}
\put(65,55){\line(0,1){2.5}}
\put(65,60){\circle{5}}
\thicklines
\put(15,10){\line(1,0){25}}
\put(15,35){\line(1,1){25}}
\put(15,10){\line(0,1){25}}
\put(40,10){\line(1,1){25}}
\end{picture}
\]
The only way to complete this to a smooth polygon $P$ is to add the
vertex $2e_1+2e_2$ indicated in the figure.  Applying
Theorem~\ref{thm:sfdegreegen} gives
\[
\deg \Sec X_P \deg \phi = 6^2 - 10\cdot 6 + 5\cdot 6 +
2\cdot 6 -12 = 6.
\]
Therefore, since $\deg \phi$ must be even, it is either 2 or 6.  If it were
6, then $\Sec X_P$ must be a linear space.  Since $\dim \Sec X_P \leq 5,$ the 
variety $\Sec X_P$ cannot fill $\PP^6.$  If $\Sec X_P$ were 
contained in a nontrivial linear space, then $X_P$ would also be contained in 
this space which is a contradiction because $X_P$ is easily seen to be nondegenerate. Therefore, we conclude that $\deg \Sec X = 2,$ and that a general point of $\Sec X_P$
lies on a unique secant line of $X_P$ by Corollary~\ref{cor:degphi}.

\medskip

\noindent{\bf Case 4.} $P$ has only edges of length 1 and dimension 3.
By Lemma~\ref{lem:dim3}, $P$ is either $\Delta_3$ or
$\Delta_1\times\Delta_2$, or else $P$ contains $P_{1,2,2}$ or
$\Delta_1^3$.  When $P$ contains $P_{1,2,2}$, we are done by the usual
combination of Lemmas~\ref{lem:projection} and~\ref{lem:goodAs}.  When
$P$ contains $\Delta_1^3$, note that Theorem~\ref{thm:3folddegreegen},
when applied to $Q = \Delta_1^3$, gives
\[
\deg \Sec X_Q \deg \phi = 6^2-21\cdot 6+48+8\cdot 8+14\cdot 8-84\cdot
0-132 = 2.
\]
This implies $\deg \phi = 2$ since $\deg \phi$ is even.  Thus a
general point of $\Sec X_Q$ lies on a unique secant line of $X_Q$.
Then $X_P$ has the same property by Lemma~\ref{lem:projection}.

\medskip

\noindent{\bf Case 5.} $P$ has only edges of length 1 and dimension
$>3$.  As usual we put a vertex of $P$ into standard position.  We
will study the 3-dimensional faces of $P$.

First suppose that $P$ has a 3-dimensional face $F$ that contains
$\Delta_1^3$ or $P_{1,2,2}$.  We can arrange for this face to lie in
$\mathrm{Cone}(e_1,e_2,e_3)$ in such a way that $e_1 + e_2 \in P$.
Now fix $i$ between $4$ and $n$ and consider the 3-dimensional face
$F_i$ of $P$ lying in $\mathrm{Cone}(e_1,e_2,e_i)$.  What are the
2-dimensional faces of $F_i$ containing $e_i$?  If both were triangles
$\Delta_2$, then the argument of the first bullet in the proof of
Lemma~\ref{lem:dim3} would imply that $F_i = \Delta_3$, which
contradicts $e_1 + e_2 \in F_i$.  It follows that one of these faces
must contain another lattice point.  This shows that there is $u_i \in
\{e_1,e_2\}$ such that $e_i + u_i \in F_i \subset P$.  If we let $B$
denote the lattice points of $\Delta_1^3$ or $P_{1,2,2}$ contained in
our original face $F$, then $P$ contains the set
\[
A = B \cup \{e_i, e_i + u_i\mid i = 4,\dots,n\}.
\]
Furthermore, the proof of Case 4 shows that a general point of $\Sec
X_B$ lies on a unique secant line of $X_B$.  By
Lemmas~\ref{lem:projection} and~\ref{lem:moregoodA}, we conclude that
a general point of $\Sec X_P$ lies on a unique secant line of $X_P$.

By Lemma~\ref{lem:dim3}, it remains to consider what happens when
every 3-dimensional face of $P$ is $\Delta_3$ or $\Delta_1 \times
\Delta_2$.  Here, Lemma~\ref{lem:3dimfaces} implies that $P =
\Delta_\ell\times\Delta_{n-\ell}$, $0 \le \ell \le n$, which cannot
occur by hypothesis.  This completes the proof of the theorem.
\end{proof}

\begin{rmk}
In Case 1 above it is possible to see that $\Sec X_A$ has the expected
dimension via results in \cite{ss}.  A lexicographic triangulation of
$A$ in which the vertices precede the other lattice points of
$\conv(A)$ will include the disjoint simplices at either end of an
edge of length $\ge 3$. This shows that $A$ satsifies the condition of
Theorem 5.4 of \cite{ss}.  (We thank S.\ Sullivant for pointing out
that a lexicographic triangulation works here.)
\end{rmk}

\subsection{Proofs of the main results} We now prove the four theorems
stated in the introduction.

\begin{proof}[Proof of Theorem~\ref{thm:dimension}]
Theorem ~\ref{thm:dimdegspecial} proves the dimensions and degrees in the
first four rows of {Table~\ref{table:dimdeg}}.  Since all of these
polytopes satisfy $\dim \Sec X_P < 2n+1$, part (1) of
Theorem~\ref{thm:dimension} follows from
Proposition~\ref{prop:degphi}.  Also, part (2) follows from
Theorem~\ref{thm:xp}.  Finally, for the last row of
{Table~\ref{table:dimdeg}}, Theorem~\ref{thm:xp} and
Corollary~\ref{cor:degphi} imply that $\dim \Sec X_P = 2n+1$ and $\deg
\phi = 2$.  Theorem~\ref{thm:degreemap} gives the desired formula for
$\deg \Sec X_P$.
\end{proof}

\begin{proof}[Proof of Corollary~\ref{cor:dimension}]
We always get the expected dimension in the last row of
{Table~\ref{table:dimdeg}}.  Since the first four rows have $\dim \Sec
X_P < 2n+1$, the only way to get the expected dimension is when $\dim
\Sec X_P = r$, where $r+1$ is the number of lattice points of $P$.
The lattice points of these polytopes are easy to count, and the
result follows.
\end{proof}

\begin{proof}[Proof of Corollaries~\ref{thm:sfdegree}
and~\ref{thm:3folddegree}] These follow from
Theorem~\ref{thm:dimension} and Corollary~\ref{cor:degphi}, together
with the formulas of Theorems~\ref{thm:sfdegreegen}
and~\ref{thm:3folddegreegen}.
\end{proof}

\subsection{Subpolytopes of $2\Delta_n$} 
We can now complete the classification of smooth subpolytopes of
$2\Delta_n$ begun in Theorem~\ref{thm:prelimclass}.

\begin{thm}
\label{thm:finalclass}
Let $P \subset 2\Delta_n$ be a smooth polytope of dimension $n$.
Then $P$ is $AGL_n(\ZZ)$-equivalent to one of
\[
\Delta_n,\ \ 2\Delta_n,\ \ (2\Delta_n)_k\ (0 \le k \le n-2),\ \
\Delta_\ell \times \Delta_{n-\ell}\ (1 \le \ell \le n-1).
\]
\end{thm}

\begin{proof} $P \subset 2\Delta_n$ gives a dominating
map $\Sec \nu_2(\PP^n) \dashrightarrow \Sec X_P$, so that $\dim \Sec
X_P \le \dim \Sec \nu_2(\PP^n) = 2n < 2n+1$.  This excludes the last
row of {Table~\ref{table:dimdeg}}, so that $P$ must be
$AGL_n(\ZZ)$-equivalent to a polytope in one of the first four rows,
as claimed.
\end{proof}

\begin{rmk}
Given a smooth lattice polytope $P \subset \RR^n$ of dimension~$n$, it
is easy to determine if $P$ is equivalent to a subpolytope of $2
\Delta_n$.  Let $d$ be the maximum number of length 2 edges incident
at any vertex of $P$ and let $v$ be any vertex that witnesses this
maximum.  Then let $Q$ be the image of $P$ under any element of
$AGL_n(\ZZ)$ that places $v$ at the origin in standard position.

We claim that $P$ is equivalent to a subpolytope of $2 \Delta_n$ if
and only if $Q$ is contained in $2 \Delta_n$.  For the nontrivial part
of the claim, note that if $P$ is equivalent to a subpolytope of $2
\Delta_n,$ then it must be equivalent to either $\Delta_{\ell} \times
\Delta_{n-\ell}$ with $1\leq \ell \leq n-1$ or $(2 \Delta_n)_k$ with
$-1 \leq k \leq n-1$ by Theorem \ref{thm:finalclass}.  The
automorphism group of $\Delta_{\ell} \times \Delta_{n-\ell}$ acts
transitively on its vertices, and the automorphism group of $(2
\Delta_n)_k$ acts transitively on vertices with a maximum number of
edges of length 2.  Therefore, by symmetry, we can check to see if $P$
is equivalent to a subpolytope of $2 \Delta_n$ at any vertex $v$ as
specified above.
\end{rmk}

\subsection{Segre-Veronese varieties}
The polytope $d_1\Delta_{n_1}\times\dots\times d_k\Delta_{n_k}$, $d_i
\ge 1$, gives an embedding of $X = \PP^{n_1} \times \cdots \times
\PP^{n_k}$ using $\mathcal{O}_X(d_1, \ldots, d_k)$.  The image $Y$ of
this embedding is a Segre-Veronese variety of dimension $n =
\sum_{i=1}^k n_i$.  We now compute the degree of $Y$ using
Theorem~\ref{thm:dimension}.

\begin{thm}
\label{thm:segver}
Let $Y$ be the Segre-Veronese of $d_1\Delta_{n_1}\times\dots\times
d_k\Delta_{n_k}$, $d_i \ge 1$.  If $\sum_{i=1}^k d_i \ge 3$,
then $\dim \Sec Y= 2n+1$ and
\begin{gather*}
\deg \Sec Y = \frac12\Big(\bigl((n_1, \ldots, n_k)!\, d_1^{n_1}\cdots
d_k^{n_k})^2 \ -\ \\ 
\sum_{\ell = 0}^n \binom{2n\!+\!1}{\ell}(-1)^{n-\ell}
\!\!\!\sum_{\sum j_i =n-
\ell}\,(n_1\!-\!j_1, \ldots,
n_k-j_k)!\, \prod_{i=1}^k \binom{n_i\! +\! j_i}{j_i}  d_i^{n_i\!-\!j_i}\Big).
\end{gather*}
where $(m_1, \ldots, m_k)! = \frac{(m_1 +\cdots + m_k)!}{m_1!\cdots
  m_k!}$ is the usual multinomial coefficient.
\end{thm}

\begin{proof} One easily sees that $d_1\Delta_{n_1}\! \times \dots
\times d_k\Delta_{n_k}$ cannot lie in $2\Delta_n$ when $\sum_{i=1}^k
d_i \ge 3$.  It follows that $\dim \Sec Y= 2n+1$ and $\deg \Sec Y$ is
given by the formula of Theorem~\ref{thm:dimension}.  It is well-known
that $\deg Y = (n_1, \ldots, n_k)!\, d_1^{n_1}\cdots d_k^{n_k}$.

To evaluate the remaining terms of the formula, we first compute
$c(T_Y)^{-1}$.  Let $H_i = c_1(\mathcal{O}_X(0, \ldots, 1, \ldots
0))$, where the 1 appears in the $i$th position.  Since $c(T_Y) =
c(T_X)$, we can use \eqref{eq:chern} to see that $c(T_Y) = \prod_{i =
1}^k (1+H_i)^{n_i+1}$.  We compute $c(T_Y)^{-1}$ by inverting each
factor and multiplying out the result.

For each $i$, we have the expansion
\begin{equation}
\label{eq:wellknown}
(1+H_i)^{-(n_i+1)}=
\sum_{j = 0}^{n_i} (-1)^j \binom{n_i+j}{j} H_i^j
\end{equation}
since $H_i^{n_i+1} = 0$.  Then we can expand the product $\prod_{i=1}^k
(1+H_i)^{-(n_i+1)}$ and take the degree $\ell$ piece to obtain
\[
c_{\ell} = (-1)^{\ell} \sum_{\sum j_i = \ell} \prod \binom{n_i
+j_i}{j_i}H_i^{j_i}.
\]

The embedding of $X$ is given by $H = d_1H_1+\cdots + d_kH_k$.  This
makes it easy to complete the computation.  For each $\ell = 0,\ldots,
n$, we need to compute
\begin{gather*}
c_{n-\ell} \cdot (d_1H_1+\cdots + d_kH_k)^{\ell} \\ =(-1)^{n-\ell}
\sum_{\sum j_i = n- \ell} \prod \binom{n_i +j_i}{j_i}H_i^{j_i} \cdot
(d_1H_1+\cdots +d_kH_k)^{\ell}.
\end{gather*}
Since the coefficient of $H_1^{n_1-j_1} \cdots H_k^{n_k-j_k}$ in
$(d_1H_1 +\cdots + d_kH_k)^{\ell}$ is $(n_1-j_1, \ldots, n_k-j_k)!\,
\prod_{i=1}^k d_i^{n_i-j_i}$, the result follows.
\end{proof}

Here are two easy corollaries of Theorem~\ref{thm:segver}.

\begin{cor}
Let $n \geq 2$. If $d \geq 3$ and $Y$ is the $d$-uple Veronese variety
of $\PP^n$, then
\[
\deg \Sec Y = \frac12\Big(d^{2n}- \sum_{j = 0}^n
(-1)^{n-j}d^j\binom{2n+1}{j} \binom{2n-j}{n-j}\Big).
\]
\end{cor}

\begin{cor}
If $n \geq 3$ and $Y$ is the $n$-fold Segre variety $\PP^1 \times
\cdots \times \PP^1$, then
\[
\deg \Sec Y = \frac12\Big((n!)^2 - \sum_{j = 0}^n
\binom{2n+1}{j}\binom{n}{n-j}j!\, (-2)^{n-j}\Big).
\]
\end{cor}

The reader should consult \cite{cgg2} for further results on the
secant varieties of Segre-Veronese varieties.

\section{Subsets of lattice points}
\label{sec:projections}

Let $P$ be a smooth polytope of dimension $n$ in $\RR^n$.  Using all
lattice points of $P$ gives the projective variety $X_P \subset
\PP^r$, $r = |P\cap\ZZ^n|-1$, in the usual way.  As in the discussion
leading up to \eqref{eq:projectionpa} in Section~\ref{sec:key}, a
subset $A \subset P \cap \ZZ^n$ gives the projective toric variety
$X_A \subset \PP^s$, $s = |A|-1$.

Here is our main result concerning the dimension and degree of the
secant variety $\Sec X_A$.  

\begin{thm}
\label{thm:xa}
Let $A \subset P \cap \ZZ^n$, where $P$ is a smooth polytope of
dimension $n$.  If $A$ contains the vertices of $P$ and their
nearest neighbors along edges of $P,$ then $\dim \Sec X_A = \dim \Sec
X_P$ and $\deg \Sec X_A$ divides $\deg \Sec X_P$.
\end{thm}

When $A$ satisfies the hypothesis of Theorem~\ref{thm:xa}, it follows
that $\dim \Sec X_A$ is given in {Table~\ref{table:dimdeg}}.
Furthermore, since $\deg \Sec X_A$ divides $\deg \Sec X_P$ and the
latter is given in {Table~\ref{table:dimdeg}}, we get an explicit
bound for $\deg \Sec X_A$.

\begin{proof}
By \eqref{eq:projectionpa} of Section~\ref{sec:key}, $A \subset
P\cap\ZZ^n$ gives a projection
\[
\pi : \PP^r \dashrightarrow \PP^s, 
\]
where $r = |P\cap\ZZ^n| - 1, s = |A| -1$.  By our hypothesis on $A$,
it is straightforward to show that $\pi$ induces an isomorphism $\pi:
X_P \to X_A$ (see, for example, the proof of the lemma on p.\ 69 of
\cite{fulton-t}).  This in turn induces a projection
\begin{equation}
\label{eq:projection}
\pi : \Sec X_P \dashrightarrow \Sec X_A
\end{equation}
since projections take lines not meeting the center to lines.  Let
$\Lambda \subset \PP^r$ be the center of the projection.  If $\Lambda$
does not meet $\Sec X_P$, then \eqref{eq:projection} is finite and
surjective, and the theorem follows easily.  Hence it remains to show
that $\Lambda \cap \Sec X_P = \emptyset$.

Suppose by contradiction that there is $z \in \Lambda \cap \Sec X_P$.
Since $A$ contains all vertices of $P$, it is easy to see that
$\Lambda$ does not meet $X_P$.  Hence $z \notin X_P$.  There are now
two cases to consider:\ either $z$ lies on a secant line, i.e., $z =
\lambda p + \mu q$, where $p \ne q \in X_P$ and $[\lambda : \mu]\in
\PP^1$, or $z$ lies on a tangent line to $p \in X_P$.

The toric variety $X_P$ is covered by open affine sets corresponding
to the vertices of $P$.  The point $p$ lies in one of these.  Assume
that we have moved this vertex into standard position at the origin so
that local coordinates near $p$ are given by sending $(t_1, \ldots,
t_n) \in \CC^n$ to $[1: t_1: t_2: \cdots : t_n: \cdots] \in \PP^r$,
and let $p = [1: \alpha_1: \alpha_2: \cdots : \alpha_n: \cdots]$.

In the first case, we may write $z = p + rq$ where 
$r = \frac{\mu}{\lambda}$ because we assume that $z\notin X.$
Since $z \in \Lambda$, all of the coordinates corresponding to the
elements of $A$ must be zero.  Therefore, $p +r q$ must be zero in
the first $n+1$ coordinates since $A$ contains each vertex and its
nearest neighbors.

Since $p$ is nonzero in the first coordinate and $p + r q$ is zero
in that coordinate, $q$ must be in the same open affine chart as $p$.
Therefore, $q = [1: \beta_1: \beta_2: \cdots: \beta_n: \cdots]$, and
then $r = -1$ and $\alpha_i = \beta_i$ for all $i$.  This implies $p
= q$, a contradiction.

Suppose now that $z$ lies on a tangent line through $p$.  From our
local coordinate system near $p$ we see that the tangent space of
$X_A$ at $p$ is spanned by $p$ and points $v_i =[0: \cdots: 0:1: 0:
\cdots: 0:*: \cdots]$.  Thus one of the first $n+1$ coordinates
is nonzero, so $z \notin \Lambda$.
\end{proof}

We can also determine when $\Sec X_A$ has the expected dimension.  

\begin{thm}\label{cor:dimensionxa}
Let $X_A \subset \PP^s$ come from $A = \{u_0,\dots,u_s\} \subset
\ZZ^n$.  If $P = \mathrm{Conv}(A)$ is smooth of dimension $n$ and $A$
contains the nearest neighbors along the edges of each vertex of $P$,
then $\Sec X_A$ has the expected dimension $\min\{s,2n+1\}$ unless $A$
is $AGL_n(\ZZ)$-equivalent to the set of all lattice points of one of
\[
2\Delta_n\ (n\ge2),\ (2\Delta_n)_k\ (0 \le k \le n-3),\
\Delta_\ell \times \Delta_{n-\ell}\ (2 \le \ell \le n-2).
\]
\end{thm}

\begin{proof}
For $P$ in the last row of {Table~\ref{table:dimdeg}}, $\dim
\Sec X_A = \dim \Sec X_P = 2n+1$ by Theorems~\ref{thm:xa}
and~\ref{thm:dimension}.  Now suppose that $P$ is
$AGL_n(\ZZ)$-equivalent to a polytope in the first four rows.  For
these polytopes, all lattice points lie on edges and all edges have
length $\le 2$.  Hence $A$ contains all lattice points of the
polytope, so Corollary~\ref{cor:dimension} applies.
\end{proof}

We conclude with two examples, one in which the degrees of $\Sec X_P$
and $\Sec X_A$ are equal and one in which they are not.

\begin{ex}
Let $P = 3 \Delta_2$ and let $A$ be the set of all lattice points in
$P$ minus the point $(1,1)$.  The variety $X_P \subset \PP^{9}$ is
the 3-uple Veronese embedding of $\PP^2$ and $X_A \subset \PP^8$ is a
projection of $X_P$ from a point.  Our results imply that $\dim \Sec
X_P = \dim \Sec X_A = 5$, $\deg \Sec X_P = 15$, and $\deg \Sec X_A$
divides 15.  A \emph{Macaulay 2} computation shows that $\deg \Sec X_A
= 15$ in this case.  Note also that a general point of $\Sec X_P$ or
$\Sec X_A$ lies on a unique secant line of $X_P$ or $X_A$,
respectively.
\end{ex}

\begin{ex}
\label{ex:hex}
Let $A = \{(0,0), (1,0), (0,1),(2,1),(1,2), (2,2)\}$ and $P$ be the
convex hull of $A$.  The set of lattice points in $P$ is $A \cup
\{(1,1)\}$.  The variety $X_P$ is a smooth toric surface of degree 6
in $\PP^6$ and $X_A$ in $\PP^5$ is the projection of $X_P$ from a
point.  Note that $X_P$ is the Del Pezzo surface obtained by blowing
up $\PP^2$ at three points, embedded by the complete anticanonical
linear system.  Our results imply that $\dim \Sec X_P = \dim \Sec X_A
= 5$ and $\deg \Sec X_P = 3$.  However, $\dim \Sec X_A = 5$ implies
that $\Sec X_A$ fills all of $\PP^5$ and hence $\deg \Sec X_A = 1$.
Thus $\pi : \Sec X_P \to \Sec X_A$ has degree 3 and a general point of
$\Sec X_A$ lies on three secant lines of $X_A$.
\end{ex}

When $X_A$ is a (possibly singular) monomial curve, the degree of
$\Sec X_A$ has been computed explicitly in the recent paper \cite{r}.

\bibliographystyle{plain}

\end{document}